%% file: LowerLDPNetworks_ArXiv.tex
\documentclass[11pt,a4paper,reqno]{amsart}

\input{packages}
\input{theorems}

\input{defs}

\usepackage{fullpage}

\title{Lower large deviations for geometric functionals} 

\begin{document}

\author{Christian Hirsch}
\author{Benedikt Jahnel}
\author{Andr{\'a}s~T{\'o}bi{\'a}s}
\address[Christian Hirsch]{University of Mannheim, Institute of Mathematics, 68161 Mannheim, Germany}
\email{hirsch@uni-mannheim.de}
\address[Benedikt Jahnel]{Weierstrass Institute for Applied Analysis and Stochastics, Mohrenstra{\ss}e 39,
10117 Berlin, Germany}
\email{jahnel@wias-berlin.de}
\address[Andr{\'a}s~T{\'o}bi{\'a}s]{Technical University of Berlin, Institute of Mathematics, Stra{\ss}e des 17.~Juni 136, 10623 Berlin, Germany.}
\email{tobias@math.tu-berlin.de}



\begin{abstract}This work develops a methodology for analyzing large-deviation lower tails  associated with geometric functionals computed on a homogeneous Poisson point process. The technique applies to characteristics expressed in terms of stabilizing score functions exhibiting suitable monotonicity properties. We apply our results to clique counts in the random geometric graph, intrinsic volumes of Poisson--Voronoi cells, as well as power-weighted edge lengths in the random geometric, $k$-nearest neighbor and relative neighborhood graph.
\end{abstract}

\maketitle

\input{result}
\input{examples}
\input{proof}

\normalsize
\input{LowerLDPNetworks_ArXiv.bbl}

\bigskip
\thanks{This work was co-funded by the German Research Foundation under Germany's Excellence Strategy MATH+: The Berlin Mathematics Research Center, EXC-2046/1 project ID: 390685689.}
%
\end{document}

%% file: packages.tex
\usepackage{tikz}
\usepackage[utf8]{inputenc}
\usepackage{amsmath,bm,bbm,amsthm, amssymb}
\usepackage{color}
\usepackage{animate}
\usepackage{etoolbox}
\usepackage{ulem}
\usepackage{comment}
\usepackage{pgf}
\usepackage{hyperref}
\usetikzlibrary{calc}
\usetikzlibrary{patterns}
\usetikzlibrary{arrows}
\usetikzlibrary{decorations.pathreplacing}
\usepackage[utf8]{inputenc}
\usepackage{dsfont}
\newcommand\numberthis{\addtocounter{equation}{1}\tag{\theequation}}

%% file: theorems.tex
\theoremstyle{plain}
\newtheorem{theorem}{Theorem}[section]

\newtheorem{lemma}[theorem]{Lemma}

\theoremstyle{definition}

\theoremstyle{remark}

\numberwithin{equation}{section}
\numberwithin{figure}{section}

%% file: defs.tex
\def\mbb{\mathbb}

\def\Z{\mbb Z}
\def\E{\mbb E}

\def\P{\mbb P}
\def\Q{\mbb Q}
\def\R{\mbb R}
\def\N{\mbb N}

\def\d{{\rm d}}

\def\a{\alpha}

\def\de{\delta}

\def\e{\varepsilon}

\def\es{\emptyset}
\def\one{\mathbbmss{1}}

\def\NN{\mathbf N}

\def\Xm{X^{-, \e}}
\def\Xmm{X^{-, M}}
\def\Xp{X^{+, \e}}
\def\Xpp{X^{+, M}}
\def\Xe{X^\e}
\def\Xee{X^M}
\def\akm{A^M}
\def\hkm{H^{\de, M}_n}
\def\ekm{E^M_n}

\def\xikm{\xi^{\de, M}}
\def\xikr{\xi^ r}

\def\vp{\varphi}
\def\kNN{k\mathrm{NN}}
\def\RN{\mathrm{RN}}
\def\ua{\uparrow}
\def\da{\downarrow}

%% file: result.tex
\section{Introduction and main results}
\label{introSec}

%
%
Considering the field of random graphs, there is a subtle difference in the understanding between upper and lower tails in a large-deviation regime. For instance, when considering the triangle count in the Erd\H{o}s--R\'enyi graph, the probability of observing atypically few triangles is described accurately via very general Poisson-approximation results \cite{janson1, janson2}. On the other hand, the probability of having too many triangles requires a substantially more specialized and refined analysis \cite{misLog}.

\begin{figure}[!htpb]
	\includegraphics[width=0.49\textwidth]{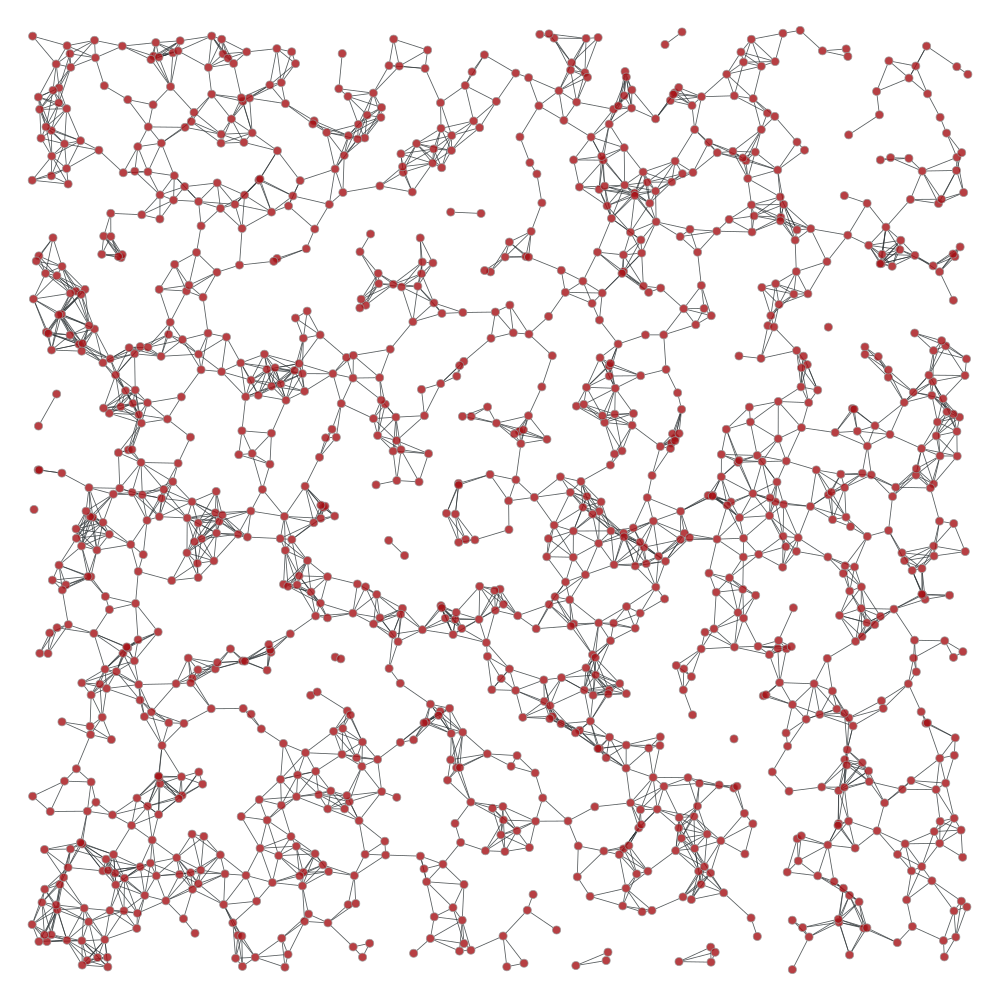}
	\includegraphics[width=0.49\textwidth]{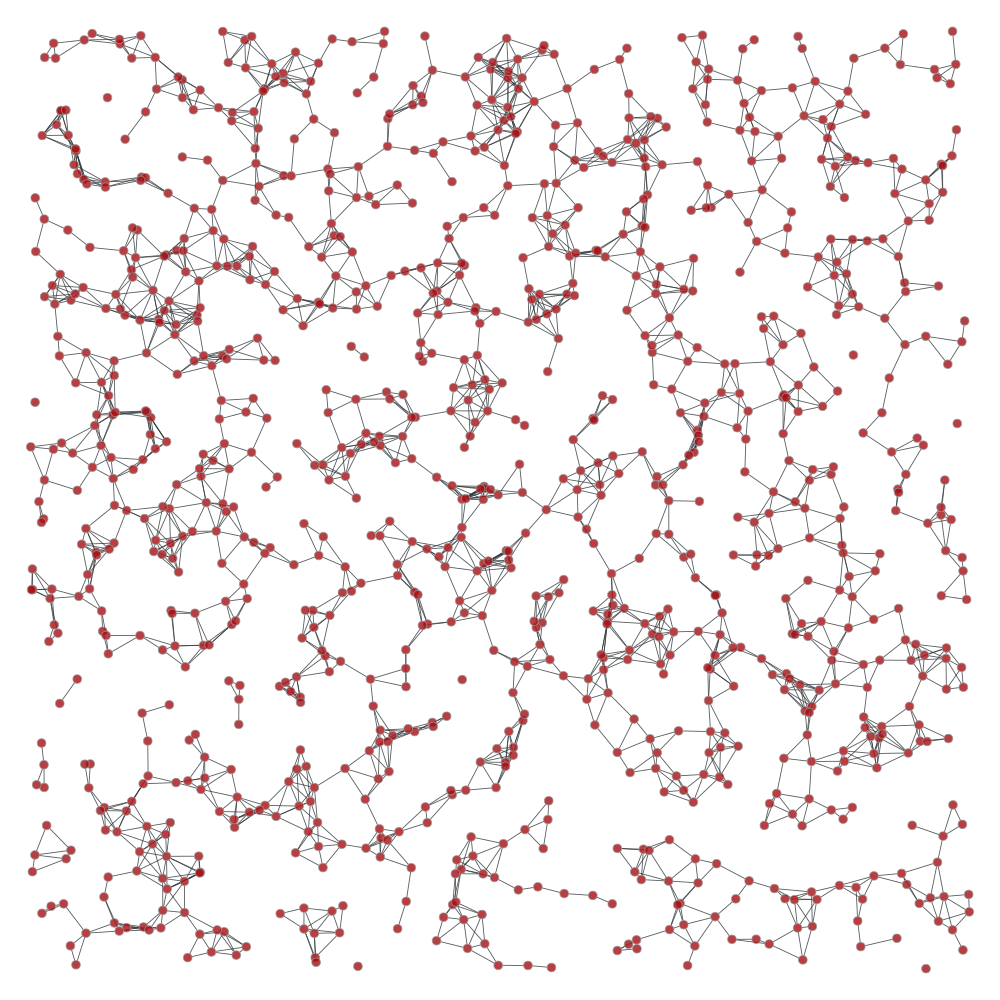}
	\caption{Typical realization of the random geometric graph (left) next to a realization having fewer than $75\%$ of the expected number of edges (right).}
	\label{rareFig}
\end{figure}

This begs the question whether a similar dichotomy also arises in the large-deviation analysis of functionals that are of geometric rather than combinatorial nature. For instance, Figure \ref{rareFig} shows a typical realization of the random geometric graph in comparison to a realization with an atypically small number of edges. In geometric probability, elaborate results are available for large and moderate deviations of geometric functionals exhibiting a similar behavior in the upper and the lower tails \cite{yukLDP2,yukLDP,eichelsbacher}. However, they prominently do not cover the edge count in the random geometric graph, whose upper tails have been understood only  recently \cite{harel}.

In the present work, we provide three general results, Theorems \ref{upThm}, \ref{gilbThm} and \ref{neighbThm}, tailored to studying large-deviation lower tails of geometric functionals. For the proofs, we resort to a method inspired by the idea of sprinkling \cite{sprinkling}. We perform small changes in those parts of the domain where the underlying point process exhibits highly pathological configurations. After this procedure, we can compare the resulting functionals to approximations that are then amenable to the point-process based large-deviation theory from \cite{georgii2} or \cite{yukLDP2,yukLDP}. Among the examples covered by our method are clique counts in the random geometric graph, inner volumes of Poisson--Voronoi cells and power-weighted edge lengths in the random geometric, $k$-nearest neighbor and relative neighborhood graph.

In the rest of this section, we set up the notation and state the main results. Then, Section \ref{sec-examples} illustrates those results through the examples. Finally, Section \ref{sec-proofs} contains the proofs.

\medskip
%
%
We study functionals on a homogeneous Poisson point process $X = \{X_i\}_{i \ge 1} \subset \R^d$ with intensity 1, whose distribution on the space $\NN$ of locally-finite configurations will be denoted by $\P$. Following the framework of \cite{yukLDP2}, these functionals are realized as averages of scores associated to the points of $X$. More precisely,  a \emph{score function} 
$$\xi:\, \R^d \times \NN \to [0, \infty)$$
is any bounded measurable function. To simplify notation, we shift the coordinate system to the considered point and write $\xi(X - X_i) = \xi(X_i, X)$. In this notation $\vp\mapsto\xi(\vp)$ acts on configurations $\vp\in \NN_o$, the family of locally-finite point configurations with a distinguished node at the origin $o \in \R^d$.

%
%
We then consider lower tails of functionals of the form
\begin{align}
	\label{hnEq}
	H_n = H_n^\xi(X) =  \frac1{n^d}\sum_{X_i \in X\cap Q_n}\xi(X - X_i),
\end{align}
i.e., averages of the score function over all points in the box $Q_n = [-n/2, n/2]^d$ of side length $n \ge 1$ centered at the origin.

%
%
In a first step, we derive upper bounds for the lower tail probabilities.  To that end, we work with approximating score functions $\xi^r$ that are  \emph{$r$-dependent} for some $r > 0$. That is, $\xi^r(\vp) = \xi^r(\vp \cap B_r)$ for every $\vp \in\NN_o$, where $B_r$ denotes the Euclidean ball of radius $r$ centered at the origin.

%
%
To state the main results, we resort to the entropy-based formulation of the large-deviation rate function. We write
$$h(\Q) =\lim_{n\ua\infty}\frac{1}{n^d}\int\d \Q_n\log\frac{\d\Q_n}{\d \P_n}$$ 
for the \emph{specific relative entropy} of a stationary point process $\Q$, where $\Q_n$ and $\P_n$ denote the restrictions of $\Q$ and $\P$ to the box $Q_n$, respectively. If $\Q_n$ is not absolutely continuous with respect to the restricted Poisson point process, we adhere to the convention that the above integral is infinite. Further, $\Q^o[\xi]$ is the expectation of $\xi$ with respect to the Palm version $\Q^o$  of $\Q$, see~\cite{georgii2} for details. Here is our first main theorem. 
%
%
\begin{theorem}[Upper bound]
	\label{upThm}
	Let $a > 0$ and assume the score function $\xi$ to be the pointwise increasing limit of a family $\{\xikr\}_{r \ge 1}$ of $r$-dependent score functions. Then, 
	\begin{align}
		\label{gilbUpEq}
		\limsup_{n \uparrow \infty} \frac {1}{n^d}\log \P(H_n \le a) \le -\inf_{\Q:\, \Q^o[\xi] \le a}h(\Q).
	\end{align}

\end{theorem}

For the lower bound, we give two sets of conditions. The first deals with score functions $\xi$ that are \emph{increasing} in the sense that $\xi(\vp) \le \xi(\psi)$ for every $\vp \subset \psi$. This applies for instance to clique counts and power-weighted edge lengths in the random geometric graph.

%
%
\begin{theorem}[Lower bound for bounded-range scores]
        \label{gilbThm}
	Let $a> 0$ and assume the score function $\xi$ to be increasing and $r$-dependent for some $r > 0$. Moreover, assume that for every $b > 0$ there exists $M = M(b) > 0$ such that $\xi(\vp) \le M$ whenever $\#\vp < b$.	Then,
	\begin{align}
		\label{gilbLowEq}
		\liminf_{n \uparrow \infty} \frac{1}{n^d}\log \P( H_n < a) \ge -\inf_{\Q:\, \Q^o[\xi] < a}h(\Q).
	\end{align}
\end{theorem}

However, many score functions are neither $r$-dependent nor increasing, or not even monotone. A prime example is the sum of power-weighted edge lengths in the $k$-nearest neighbor graph, see Section~\ref{sec-examples}. 
Still, this example and many other score functions are stabilizing, $R$-bounded and weakly decreasing in the following sense. 

%
%
First, a score function $\xi$ is \emph{stabilizing} if there exists a $\P^o$-almost surely finite measurable \emph{stabilization radius} $R:\ \NN_o \to [0, \infty]$, such that $\{R(X) \le r\}$ is measurable with respect to $X \cap B_r$ for every $r \ge 0$ and
$$\P^o\big( \xi(X) = \xi(X \cap B_{R(X)})\big) = 1.$$
In words, $\xi(X)$ does not depend on the configuration outside the ball $B_{R(X)}$. We call $R$ {\em decreasing} if $R(\vp\cup\{x\})\le R(\vp)$ for all $\vp\in \NN_o$ and $x\in \R^d$. 

%
%
Second, $\xi$ is \emph{$R$-bounded} if for every $\de > 0$ and sufficiently large $M=M(\delta) \ge 1$,
$$\P^o\big(\{R(X) \le M\} \cap \{\xi(X) \ge \de M^d\}\big) = 0.$$
Loosely speaking, the score function is negligible compared to the $d$th power of the stabilization radius.

%
%
Third, $\xi$ is \emph{weakly decreasing} if 
$$
\P\big(\#\{y\in X\colon \xi(X\cup\{o\} - y)>\xi(X - y)\}\le k\big) = 1
$$
holds for some $k \ge 1$.
In words, for all but at most $k$ points of a configuration, adding a new point to the configuration decreases the score function value of the point.

Finally, we need to ensure that sprinkling a sparse configuration of Poisson points yields control on the stabilization radii of the points in a box. More precisely, we assume that the stabilization radius is \emph{regular} in the following sense.
Let $\Xpp$ denote a Poisson point process with intensity $M^{-d}$ that is independent of $X$.
 Then, we assume that there exists $K_0 > 0$ with the following property. For every $\de > 0$ there exist $M_0 = M_0(\de) \ge 1$ and $n_0=n_0(\de) \geq 1$ such that for all $M  \ge M_0$ and $n \geq n_0$,
$$\P\big(\{ \Xpp(Q_n) \le K_0 (n/M)^d \} \cap  E_n^{M, +}|X  \big) \ge \exp(-\de n^d)$$
holds almost surely. Here, 
for $\vp\in\NN$ and any measurable subset $A\subset \R^d$, we write $\vp(A)=\#\{x\in \vp\colon x\in A\}$ for the number of points of $\vp$ contained in $A$, and
$$E_n^{M, +} = \max_{X_i \in (X \cup \Xpp)\cap  Q_n} R\big((X \cup \Xpp) - X_i\big) \le M$$
denotes the event that after the sprinkling, the stabilization radii of all points in $Q_n$ are at most $M$. Here is the corresponding main result.

%
%
\begin{theorem}[Lower bound for stabilizing scores]
	\label{neighbThm}
	Let $a > 0$ and $\xi$ be a weakly-decreasing $R$-bounded score function with a decreasing and regular radius of stabilization. Then, \eqref{gilbLowEq} remains true.
\end{theorem}

%% file: examples.tex
\section{Examples}\label{sec-examples}
In this section, we discuss how to apply the results announced in Section \ref{introSec} to a variety of examples arising in geometric probability. More precisely, Sections \ref{rggSec}, \ref{vorSec} and \ref{knnSec} are devoted to characteristics for the random geometric graph, the Voronoi tessellation, $k$-nearest neighbor graphs and relative neighborhood graphs, respectively.

%
%
\subsection{Clique counts and power-weighted edge lengths in random geometric graphs} 
\label{rggSec}
As a first simple application of our results, consider the set 
$$
C_k(\vp) = C_{k, t}(\vp)= \big\{\{x_1,\dots, x_k\}\subset\vp\colon x_1=o \text{ and }|x_i - x_j|<t \text{ for all $i \ne j$}  \big\}
$$
of \emph{$k$-cliques} associated to the origin in the geometric graph on $\vp \in \NN_o$ with connectivity radius $t > 0$. 
Then, for $k \ge 2$ and $\a \ge 0$, the score functions
$$
\xi_k(\vp) = \frac 1k \#C_k(\vp) \qquad\text{ and }\qquad \xi'_\a(\vp) = \frac 12\sum_{\substack{x \in \vp\colon  |x| < t}}|x|^\a
$$
count the number of $k$-cliques containing the origin and the power-weighted edge lengths at the origin, respectively.
Note that $\xi_k$ and $\xi'_\a$ are $t$-dependent and increasing. Additionally, if $\#\vp < b$, then $\xi_k(\vp) \le k^{-1}b^{k - 1}$ and $\xi'_\a(\vp) \le  t^\a b$. Hence Theorems \ref{upThm} and \ref{gilbThm} are applicable. 

%
%
Further examples arise in the context of topological data analysis. More precisely, the number of $k$-cliques containing the origin is precisely the number of $k$-simplices of the Vietoris--Rips complex containing the origin. Similar arguments also apply to the \v{C}ech complex, the second central simplicial complex in topological data analysis. We refer the reader to \cite[Section 2.5]{chazal} for precise definitions and further properties.

%
%
\subsection{Intrinsic volumes of Voronoi cells} 
\label{vorSec}
Recall the definition of the Voronoi cell at the origin of a locally-finite configuration $\vp\in\NN_o$, i.e., 
$$C_o(\vp) = \{x \in \R^d:\,|x| \le \inf_{y \in \vp}|x - y|\}.$$
Recall that since $C_o(\vp)$ is a convex body, its \emph{intrinsic volumes} $v_0(C_o), v_1(C_o), \dots, v_d(C_o)$ can be computed. They are key characteristics of a convex set, e.g., $v_1$, $v_{d - 1}$ and $v_d$ are proportional to the mean width, the surface area and the volume, respectively. We refer the reader to \cite[Section 14.2]{sWeil} for a precise definition and further properties. In particular, considering $v_1$ in dimension $d = 2$, the associated characteristic $n^dH_n$ becomes the total edge length of the Voronoi graph, so that we obtain a link to the setting studied in \cite[Section 2.4.1]{yukLDP}. Due to the intricate geometry, deriving a full large deviation principle even for a strictly concave function of the edge length was only achieved for a Poisson point process that is restricted to a lattice instead of living in the entire Euclidean space. This example illustrates that even in situations where understanding the large-deviation upper tails requires a delicate geometric analysis, the lower tails may be more accessible.

%
%
More precisely, consider the score functions
$$
\xi_k(\vp) = v_k(C_o(\vp))
$$
and note that $\xi_k^r(\vp) = v_k\big(C_o(\vp)\cap B_r\big)$ is a $4r$-dependent, pointwise increasing approximation of $\xi_k(\vp)$. Hence, the upper bound of Theorem \ref{upThm} applies.

%
%
For the lower bound, the conditions of Theorem \ref{neighbThm} can be satisfied using the following definitions. The radius of stabilization is described in \cite[Section 6.3]{gaussLim}: Take any collection $\{S_i\}_{i\in I}$ of cones with apex at the origin and angular radius $\pi/12$ whose union covers $\R^d$, where $I = I(d) \in \N$. Let $S_i^+$ denote the cone that has the same apex and symmetry hyperplane as $S_i$ and has the larger angular radius $\pi/6$. 
Then, we define the stabilization radius 
\begin{align}\label{StabRadVor}
R(\vp) = 2\max_{i\in I}\min_{x\in \vp\cap S_i^+}|x|,
\end{align}
as twice the radius at which the origin has a neighbor in every extended cone. In particular, both $R$ and $\xi_k$ are decreasing. Since $C_o(\vp) \subset B_{R(\vp)}$, we deduce that 
$$\xi_k(\vp) \le v_k(B_{R(\vp)}) = R(\vp)^k v_k(B_1).$$
In particular, $\xi_k$ is $R$-bounded for $k < d$. Finally, we define for a suitable constant $L = L(d) \ge 1$ the event 
\begin{equation}
	\akm_n=\{\Xpp(Q_{M/L}(z)) = 1\text{ for all }z\in (M/L)\Z^d\cap Q_{2n}\}
\end{equation}
that $\Xpp$ has precisely one point in each sub-box from an $M/L$-partition of the box $Q_{2n}$. It follows from the definition of $R$ that the event $E_n^{M, +}$ occurs whenever $\akm_n$ occurs, provided that $L$ is chosen sufficiently large. 
Moreover, setting $K_0 = (2L)^d$, we deduce that $\Xpp(Q_n) \le K_0 (n/M)^d$ under $\akm_n$. Hence, it remains to establish the asserted lower bound on the probability $\P(\akm_n)$. Fixing $\de > 0$ and invoking the independence property of the Poisson point process yields that 
$$\P(\akm_n) = \P(\Xpp(Q_{M/L}) = 1)^{(2nL/M)^d} = {\rm e}^{-(2n/M)^d}L^{-(2nL/M)^d} \ge {\rm e}^{-\de n^d},$$
provided that $M = M(\de)$ is sufficiently large. Summarizing the above findings, we deduce that Theorem \ref{neighbThm} can be applied to get the lower bound on the rate function.

%
%
\subsection{Power-weighted edge counts in $k$-nearest neighbor graphs and relative neighborhood graphs}\label{knnSec}
Finally, we elucidate how to apply Theorem \ref{neighbThm} to the power-weighted edge count of two central graphs in computational geometry, namely the $k$-nearest neighbor graph and the relative neighborhood graph. As we shall see, in contrast to the Voronoi example presented in Section \ref{vorSec}, we encounter here score functions that are weakly decreasing but not decreasing. A full large deviation principle for the total edge length of the $k$-nearest neighbor graph is described in \cite[Section 2.3]{yukLDP}, and we believe that the proof should extend to power-weighted edge lengths with a power strictly less than $d$. Nevertheless, we apply here our approach towards the large-deviation lower tails as it can be directly adapted to the bidirectional $k$-nearest neighbor graph, the relative neighborhood graph and possibly further graphs.

%
%
In the \emph{undirected $k$-nearest neighbor graph}, $\xi$ expresses the powers of distances between any point and the origin, such that at least one of them belongs to the set of $k$ nearest neighbors of the other one. To be more precise, 
\begin{align}\label{k-radius}
	\mathfrak R_k(\vp) =\inf \{ r>0 \colon \vp(B_r)  \geq k + 1 \}
\end{align}
defines the \emph{$k$-nearest neighbor radius} of $o$ in $\vp \in \NN_o$.
Then, for some $\a \ge 0$, the score function corresponding to the sum of power-weighted edge lengths of the $k$-nearest neighbor graph is defined via 
\begin{align*}
	\xi_{k, \a}(\vp) = \frac12 \sum_{\substack{x \in \vp\colon |x| \le \mathfrak R_k(\vp) \vee \mathfrak R_k(\vp-x)}}|x|^\a.
\end{align*}
In particular, we recover the number of edges by setting $\a = 0$. As noted in \cite[Section 6.3]{gaussLim}, to construct a radius of stabilization we can proceed as in \eqref{StabRadVor} except for replacing $\min_{x\in \vp\cap S^+_i}|x|$ by the distance of the $k$th closest point from the origin in $\vp \cap S_i^+$. Hence, $\xi_{k, \a}$ becomes stabilizing with a decreasing stabilization radius. In the same vein, a minor adaptation of the arguments in Section \ref{vorSec} yield the regularity and $R$-boundedness for $\a < d$.

In order to apply Theorem~\ref{neighbThm} for the lower bound, it remains to verify the following.
\begin{lemma}
	\label{lemma-wd_undirknng}
$\xi_{k, \a}$ is weakly decreasing.
\end{lemma}
\begin{proof}
Let us call $\vp\in\NN$ \emph{nonequidistant} if for all $y,z,v,w \in \vp$, $|y-z| = |v-w| > 0$ implies $\{ y,z \}=\{ v, w \}$. 
First note that for any $x\in\R^d$, under $\P$, almost all configurations $\vp\cup\{x\}$ are nonequidistant. 
We claim that for any nonequidistant configuration $\vp\cup\{x\}$, we have
for all but at most $k$ points $y \in \vp$ that
\[ \xi_k(\vp \cup \{ x \} -y) \leq \xi_k(\vp - y). \numberthis\label{NNdecrease} \]
 Indeed, for $y \in \vp$, let us define the set of $k$ nearest neighbors of $y$ in $\vp$ as follows
	\[ \kNN(\vp,y)=\big(B_{\mathfrak R_k(\vp-y)}(y) \cap \vp\big) \setminus \{ y \}. \]
Now,
if $y \in \kNN(\vp\cup \{ x \},x)$, then possibly $\xi_k(\vp \cup \{ x \} -y) > \xi_k(\vp - y)$. We claim that else \eqref{NNdecrease} holds. Indeed, if $y \notin \kNN(\vp\cup \{ x \},x)$, then there are two possibilities. If $x \in \kNN(\vp\cup \{ x \},y)$, then $x$ replaced precisely one neighbor $z$ of $y$ and is closer to $y$ than $z$. More precisely, note that $|x-y|\leq \mathfrak R_k(\vp\cup\{x\}-y) \leq \mathfrak R_k(\vp-y)$. Hence, there exists $z \in \kNN(\vp,y)$ such that $|z-y|=\mathfrak R_k(\vp-y)$ and
$z \notin  \kNN(\vp \cup \{ x \},y)$, the neighbor of $y$ that is replaced by $x$. Additionally, for any $w \in \kNN(\vp,y) \setminus \{ z \}$ also $w \in \kNN(\vp \cup \{ x \},y)$. Further, also for any $v\in\vp$ such that $y\in\kNN(\vp\cup\{x\},v)$ we have $y\in\kNN(\vp,v)$. Hence,
\[ \xi_k(\vp \cup \{ x \} -y ) - \xi_k(\vp-y)\le  |x-y|^{\a} - |z-y|^{\a} \leq 0, \]
which is \eqref{NNdecrease}.
The other possibility is that $x \notin \kNN(\vp\cup \{ x \},y)$. Then the addition of $x$ can only remove edges that were present due to the fact that some other point had $y$ as a neighbor. In this case, $\xi(\vp \cup \{ x \}-y)=\xi(\vp-y)$ unless there exists $z \in \vp$ such that $y \in \kNN(\vp,z)$ but $y \notin \kNN(\vp \cup \{ x \},z)$, which must be due to the property that $x \in \kNN(\vp \cup \{x\},z)$. So again, the addition of $x$ can only remove such an edge and hence again~\eqref{NNdecrease} holds for $y$. 
\end{proof}
Note that the approach presented above also applies to further graphs studied in computational geometry. The most immediate adaptation concerns the \emph{bidirectional $k$-nearest neighbor graph}, see~\cite{BB08}, where in the definition of the score function, we replace $\mathfrak R_k(\vp) \vee \mathfrak R_k(\vp - x)$ by $\mathfrak R_k(\vp) \wedge \mathfrak R_k(\vp - x)$. Not only can we take the same radius of stabilization, but also Lemma~\ref{lemma-wd_undirknng} remains valid. As a third example, we showcase the \emph{relative neighborhood graph}. Here, for $\a\ge 0$ and $\vp\in\NN_o$ the score function is given by
\begin{align*}
	\xi_{\rm RN}(\vp) = \frac12 \sum_{\substack{x \in \vp\colon \vp\cap B_{|x|}(o)\cap B_{|x|}(x)=\emptyset}}|x|^{\a}.
\end{align*}
The {\em relative neighborhood graph} is a sub-graph of the Delaunay tessellation, and in fact we can reuse the radius of stabilization from Section \ref{vorSec}. Finally, proving the analog of Lemma \ref{lemma-wd_undirknng} reduces to the observation that the degree of every node in the relative neighborhood graph is bounded by a constant $K=K(d)$, see \cite[Section IV]{RNG}.
What remains to be verified is that $\xi_{\rm RN}$ is weakly decreasing. 
\begin{lemma}\label{lemma-almostincreasing_RN}
$\xi_{\rm RN}$ is weakly decreasing.
\end{lemma}
\begin{proof}
We claim that for any nonequidistant configuration $\vp\cup\{x\}$ with $\vp\in\NN$, 
for all but at most $K$ points $y \in \vp$,
\[ \xi_{\rm RN}(\vp \cup \{ x \} -y) \leq \xi_{\rm RN}(\vp - y)\numberthis\label{RNdecrease} \]
holds. Indeed, for $y \in \vp$, let us define the set of relative neighbors of $y$ in $\vp$ as follows
\[ \RN(\vp,y):=\{z\in \vp\setminus\{y\}\colon \vp\cap B_{|z-y|}(y)\cap B_{|z-y|}(z)=\emptyset\}, \]
and note that $z\in \RN(\vp,y)$ if and only if $y\in\RN(\vp,z)$. In particular, $\#\RN(\vp,y)\le K$ for any $y\in \vp$. So, if $y \in \RN(\vp\cup \{ x \},x)$, then possibly $\xi_{\rm RN}(\vp \cup \{ x \} -y) > \xi_{\rm RN}(\vp - y)$. But if $y \notin \RN(\vp\cup \{ x \},x)$, then 
\begin{align*}
&\xi_{\rm RN}(\vp \cup \{ x \} -y ) - \xi_{\rm RN}(\vp-y)\\
&=\tfrac{1}{2}\sum_{z\in \vp-y}|z-y|^\a\Big(\one\{z\in \RN(\vp\cup\{x\},y)\}-\one\{z\in \RN(\vp,y)\}\Big)\le0,
\end{align*}
as asserted.
\end{proof}

%% file: proof.tex
\section{Proofs}\label{sec-proofs}

In this section we provide the proofs of the main theorems. 

%
%
\subsection{Proof of Theorem \ref{upThm}}
\label{upSec}
The proof of the upper bound relies on the level-3 large deviation principle for the Poisson point process from \cite[Theorem 3.1]{georgii2}.
%
%
\begin{proof}[Proof of Theorem \ref{upThm}]
	Replacing $\xikr$ by $\xikr \wedge r$ if necessary, we may assume that $\xikr$ is bounded above by $r$. Then, $\xikr$ is a bounded local observable, so that by the contraction principle \cite[Theorem 4.2.10]{dz98}  and \cite[Theorem 3.1]{georgii2},
	$$\limsup_{n \ua \infty} \frac1{n^d} \log \P(H_n \le a) \le \limsup_{n \ua \infty} \frac1{n^d} \log \P(H_n^{\xikr} \le a) \le -\inf_{\Q:\, \Q^o[\xikr] \le a}h(\Q).$$
	Hence, it suffices to show that 
	$$-\lim_{r \ua \infty}\inf_{\Q:\, \Q^o[\xikr] \le a}h(\Q) \le -\inf_{\Q:\, \Q^o[\xi] \le a}h(\Q).$$
	Let $\{\Q_k\}_{k \ge 1}$ be a family of stationary point processes  such that $\Q_k^o[\xi^k] \le a$ and 
	$$\lim_{k \ua \infty}h(\Q_k) = \lim_{r \ua \infty}\inf_{\Q:\, \Q^o[\xikr] \le a}h(\Q).$$
	Let $\Q_*$ be a subsequential limit of $\{\Q_k\}_{k \ge 1}$. To simplify the presentation, we may assume  $\Q_*$ to be the limit of $\{\Q_k\}_{k \ge 1}$. Then, by monotone convergence,
		$$\Q_*^o[\xi] \le  \lim_{r \ua \infty}\Q_*^o[\xikr] = \lim_{r \ua \infty}\lim_{k \ua \infty}\Q_k^o[\xi^r]\le \limsup_{k \ua \infty}\Q_k^o[\xi^k] \le a.$$
Since the specific relative entropy $h$ is lower semicontinuous, we arrive at 
		$$\liminf_{k \ua \infty} h(\Q_k) \ge h(\Q_*) \ge \inf_{\Q:\, \Q^o[\xi] \le a}h(\Q),$$
		as asserted.
\end{proof}

\subsection{Proof of Theorem~\ref{gilbThm}}

To prove  Theorem \ref{gilbThm}, we consider the truncation $\xi^M=\xi \wedge M$ of the original increasing and $r$-dependent score function $\xi$ at a large threshold $M > 1$ and write $H_n^M=H_n^{\xi^M}$.
%
%
In comparison to the arguments in Section \ref{upSec}, the proof of the lower bound is more involved, since we can no longer replace $\P(H_n \le a)$ by $\P(H_n^M \le a)$.
Instead, we rely on a sprinkling approach. For this method to work, we need that the total number of points in pathological areas is small with high probability. More precisely, we say that a point $X_i \in X$ is \emph{$b$-dense} if $X(Q_r(X_i)) > b$ and write 
	$$N_{b, n}=N_{b,n}(X) = \#\{X_i \in X\cap Q_n:\, \text{$X_i$ is $b$-dense}\}$$
	for the total number of $b$-dense points in $Q_n$. Then, $b$-dense points are indeed rare.
\begin{lemma}[Rareness of $b$-dense points]
	\label{rareDenseLem}
	Let $\de > 0$. Then,
$$		\limsup_{b \ua \infty}\limsup_{n \ua \infty}\frac1{n^d} \log \P(N_{b, n} > \de n^d) = -\infty.$$
\end{lemma}

%
%
In the second step, we remove all $b$-dense points through the coupling. That is, we let $\Xm$ be an independent thinning of $X$ with survival probability $1 - \e$. Furthermore, we let $\Xp$ be an independent Poisson point process with intensity $\e > 0$. Then, the coupled process 
	$$\Xe = \Xm \cup \Xp$$
	is again a Poisson point process with intensity 1.
	Now, let
	$$E_{b, n} = \{\Xp \cap Q_n = \es\} \cap \{\Xm \cap Q_n \text{ has no $b$-dense points}\}$$
	be the event that $\Xp$ has no points in $Q_n$ and that $\Xm$ does not contain any $b$-dense points in $Q_n$.
\begin{lemma}[Removal of $b$-dense points]
	\label{remDenseLem}
	Let $b, n, \e > 0$. Then, $\P$-almost surely, 
	$$\P(E_{b,n}|X) \ge  \exp(-\e n^d + {N_{b, n}}\log(\e)).$$
\end{lemma}

	%
	%
Before showing Lemmas \ref{rareDenseLem} and \ref{remDenseLem}, we illustrate how they enter the proof of \eqref{gilbLowEq}.
\begin{proof}[Proof of Theorem \ref{gilbThm}]
	Let $M > 0$. Then, by \cite[Theorem 3.1]{georgii2},
	$$
	\liminf_{n \ua \infty} \frac1{n^d} \log \P(H_n^M< a) \ge  -\inf_{\Q:\, \Q^o[\xi^M] < a }h(\Q) \ge -\inf_{\Q:\, \Q^o[\xi] < a }h(\Q).
	$$
	Hence, it remains to show that 
	\begin{align}
	\label{gilbLowEq2}
	\liminf_{n \ua \infty} \frac1{n^d} \log \P(H_n < a) \ge \liminf_{M \ua \infty}\liminf_{n \ua \infty} \frac1{n^d} \log \P(H_n^M < a). 
	\end{align}
	
	%
	%
	
	%
	%
	Let $b, \de, \e > 0$ be arbitrary. Now, since $\xi$ is increasing,
	$$\P(H_n < a) = \P(H_n(\Xe) < a) \ge \P(\{H_n^{M(b)} < a\} \cap E_{b, n}) = \E\big[\one\{H_n^{M(b)} < a \}\P[E_{b,n}\,|\,X]\big].$$
	Thus, by Lemma \ref{remDenseLem}, 
	\begin{align*}
		\P(H_n < a)
		&\ge \exp(-\e n^d)\E\big[\one\{H_n^{M(b)} < a \}\e^{N_{b, n}}\big]\\
		&\ge \exp\big((\de \log(\e) - \e) n^d\big)\P(H_n^{M(b)} < a) - \P(N_{b, n} > \de n^d).
	\end{align*}
	Since $X$ and $\Xe$ share the same distribution, Lemma \ref{rareDenseLem} allows us to choose $b = b(\de) > 0$ sufficiently large such that 
	$$\liminf_{n \ua \infty} \frac1{n^d}\log\P(H_n < a) \ge \de \log(\e) - \e + \liminf_{n \ua \infty}\frac1{n^d} \log\P(H_n^{M(b)} < a).$$
	Hence, sending $\e\da0$, $\de \da 0$,  and $b \ua \infty$ concludes the proof of \eqref{gilbLowEq2}.
\end{proof}

%
%
\begin{proof}[Proof of Lemma \ref{rareDenseLem}]
	Consider a subdivision of $Q_n$, for sufficiently large $n \ge 1$, into sub-boxes $Q_a(z_i)=z_i+Q_a$ of side length $a>r$ where $z_i\in a\Z^d$. Let $N_i = X(Q_a(z_i))$ be the number of points in the $i$th sub-box and $N'_i =X(Q_{3a}(z_i))$ be the number of points the $i$th sub-box plus its adjacent sub-boxes. Then, $N_{b,n}\le N_{b,n}''$, where
	$$N_{b,n}'' = \sum_{i\in a\Z^d\cap Q_n}N_i\one\{N'_i>b\},$$
so that by the exponential Markov inequality, for all $t> 0$, 
\begin{align*}
	\log \P(N_{b, n} > \de n^d)&\le \log \P(N_{b, n}'' > \de n^d)\le  -\de t n^d+\log \E[\exp(tN_{b, n}'')].
\end{align*}
Since the random variables $N_i\one\{N'_i>b\}$ and $N_j\one\{N'_j>b\}$ are independent whenever $\Vert z_i - z_j\Vert_\infty\ge 3$, we have $3^d$ regular sub-grids of $a\Z^d$ containing independent random variables $N_i\one\{N'_i>b\}$. Thus, using H\" older's inequality, independence and the dominated convergence theorem, we arrive at 
\begin{align*}
	\limsup_{b \ua \infty}\limsup_{n \ua \infty}\frac{1}{n^d}\log \E[\exp(tN_{b,n}'')]\le\frac{1}{(3a)^{d}}\limsup_{b \ua \infty} \log \E\big[\exp(3^dtN_o\one\{N'_o>b\})\big]
	 =\frac{1}{(3a)^{d}}.
\end{align*}
Since $t > 0$ was arbitrary, we conclude the proof.
\end{proof}

%
%
\begin{proof}[Proof of Lemma \ref{remDenseLem}]
	First, since $\Xp$ and $\Xm$ are independent, it suffices to compute 
	$$\P(\Xp \cap Q_n = \es\,|\,X)\qquad \text{ and }\qquad \P(\Xm \cap Q_n \text{ has no $b$-dense points}\,|\,X)$$ separately. The void probabilities for a Poisson point process give that $$\P(\Xp \cap Q_n = \es\,|\,X) = \exp(- \e n^d).$$ 
	Next, since $\Xm$ is an independent thinning of $X$ with probability $\e$, we arrive at 
	$$\P(\Xm \cap Q_n \text{ has no $b$-dense points}\,|\,X) \ge \e^{N_{b, n}},$$
	as asserted.
\end{proof}

\subsection{Proof of Theorem \ref{neighbThm}}\label{sec-stabproof}

In order to prove the lower bound for stabilizing score functions, we use sprinkling to regularize sub-regions that are not sufficiently stabilized.
%
%
	%
	%
	Let us define the approximation 
	$$\xikm(\varphi) = \xi(\varphi \cap Q_M) \wedge \de M^d$$
	and write $\hkm = H_n^{\xikm}$.

	%
	%
	Similarly as before, we consider a coupling construction. Now, we let $\Xmm$ denote an independent thinning of $X$ with survival probability $1 - M^{-d}$ and $\Xpp$ an independent Poisson point process with intensity $M^{-d}$. Then, 
	$$\Xee = \Xmm \cup \Xpp$$ 
	defines a unit-intensity Poisson point process.

	%
	%
	In this coupling, we consider events in which the sprinkling $\Xpp$ adds points wherever necessary to reduce the stabilization radius. More precisely, let
	$$\ekm = \{\Xmm\cap Q_n = X\cap Q_n\} \cap \big\{\Xpp(Q_n) \le K_0 (n/M)^d\big\} \cap E^{M, +}_n.$$
	As we shall prove below, the events $\ekm$ occur with a high probability.

	%
	%
	\begin{lemma}[Sprinkling regularizes with high probability]
		\label{remreachinglem}
Let $\de > 0$ and $n\ge M\ge 1$ sufficiently large. Then, under the assumptions of Theorem \ref{neighbThm}, $\P$-almost surely, 
		\[\P(\ekm| X) \ge \exp\big(X(Q_n)\log(1-M^{-d}) - \de n^d\big). \]
\end{lemma}
\begin{proof}
	Indeed, for given $X$, the event $\{\Xmm\cap Q_n=X\cap Q_n\}$ has probability $(1-M^{-d})^{X(Q_n)}$ and is independent of the event $\big\{\Xpp(Q_n) \le K_0  (n/M)^d\big\} \cap E^{M, +}_n$, which has probability at least $\exp(-\de n^d)$. 
\end{proof}

Now, we conclude the proof of Theorem~\ref{neighbThm}.
\begin{proof}[Proof of Theorem~\ref{neighbThm}]
	Let $\de>0$ and $M=M(\delta) > 1$ sufficiently large. Then, by $R$-boundedness,
	$$\P(H_n < a) = \P(H_n(\Xee) < a) \ge \P\big(\{\hkm(\Xee) < a\} \cap \ekm\big).$$
	Moreover, under the event $\ekm$,
	\begin{align*}
		\hkm(\Xee) &=\frac1{n^d}\sum_{X_i \in \Xpp\cap Q_n} \xikm(\Xee - X_i)  +\frac1{n^d}\sum_{X_i \in X\cap Q_n} \xikm(\Xee - X_i) \\
		&\le K_0\de +\hkm(X) +   \frac1{n^d}\sum_{X_i \in X \cap Q_n} \big(\xikm(\Xee - X_i)  -\xikm(X - X_i)\big).
	\end{align*}
Let us write $X^{M,0} = X$ and $X^{M,j+1} = X^{M,j}\cup\{\Xpp_j\}$ where $\{\Xpp_j\}_{1 \le j \le N(M)}$ is an arbitrary ordering of $\Xpp$. Then, since $\xi$ is weakly decreasing,
\begin{align*}
&  \sum_{X_i \in X \cap Q_n} \big(\xikm(\Xee - X_i)  -\xikm(X-X_i)\big) \\
&\quad= \sum_{X_i \in X \cap Q_n}\sum_{j \le N(M)} (\xikm(X^{M,j} - X_i)  -\xikm(X^{M,j-1}-X_i))\\
&\quad\le\de M^d\sum_{j \le N(M)}\sum_{X_i \in X \cap Q_n}\one\big\{\xikm(X^{M,j} - X_i) >\xikm(X^{M,j-1}-X_i)\big\}\\
&\quad\le k \de M^dN(M).
 \end{align*} 
	Further note that $N(M) \le K_0(n/M)^d$, and thus we arrive at
	$$\P(H_n(\Xee) < a) \ge \P\big(\{\hkm(\Xee) < a\} \cap \ekm\big) \ge \P\big(\{\hkm(X) < a - 2kK_0 \de\} \cap \ekm\big).$$
	Now, by conditioning on $X$ and applying Lemma~\ref{remreachinglem} for sufficiently large $n\ge M\ge 1$, 
	\begin{align*}
		&\P(H_n(\Xee) < a) \ge \E\big[\one\{\hkm(X) < a - 2kK_0 \de \} \P(\ekm\, |\, X)\big]\\
		&\quad\ge\E\big[\one\{\hkm(X) < a - 2kK_0 \de\}\exp\big(X(Q_n)\log(1-M^{-d})\big)\big] \exp(-\de n^d).
	\end{align*}
	Moreover, for any $c > 0$, 
		\begin{align*}
			&\E\big[\one\{\hkm(X) < a - 2kK_0 \de\}\exp\big(X(Q_n)\log(1-M^{-d})\big)\big]\\
		&\quad\ge\exp\big(cn^d\log(1-M^{-d})\big)\P(\{\hkm(X) < a - 2kK_0 \de\}\cap\{X(Q_n)<cn^d\}), 
	\end{align*}
	where for the first factor, 
			\begin{align*}
				\liminf_{M\ua\infty}\frac1{n^d}\log\big(\exp\big(cn^d\log(1-M^{-d})\big)\big) = \liminf_{M\ua\infty}c\log(1-M^{-d}) = 0.
	\end{align*}
	Now, for the second factor, 
		\begin{align*}
			\P\big(\{\hkm(X) < a - 2kK_0 \de\}\cap\{X(Q_n)<cn^d\}\big)		&\ge \P(\hkm(X) < a - 2kK_0 \de) \\
			&\phantom= -\P(X(Q_n)\ge cn^d),
	\end{align*}	
	where for large $c$ the second summand plays no role in the large deviations. Applying \cite[Theorem 3.1]{georgii2} on the local bounded observable $\xikm$ yields that
	\begin{align*}
		\liminf_{n \ua \infty}\frac1{n^d}\log\P\big(\hkm(X) < a -  2kK_0 \de\big) \ge -\inf_{\Q:\, \Q^o[\xikm] < a -  2kK_0 \de} h(\Q).
	\end{align*}
	Finally, if $\Q^o[\xi] < a$, then $\limsup_{M \ua \infty}\Q^o[\xi^{\delta,M}] < a - 2kK_0 \de$ for a sufficiently small $\delta>0$, so that
	$$\liminf_{M \ua \infty} \Big( -\inf_{\Q:\, \Q^o[\xi^{\delta,M}] < a -  2kK_0 \de} h(\Q) \Big) \ge -\inf_{\Q:\, \Q^o[\xi] < a}h(\Q),$$
	 as asserted.
\end{proof}

%% file: LowerLDPNetworks_ArXiv.bbl
\normalsize
\newcommand{\etalchar}[1]{$^{#1}$}